\documentclass[10pt,a4paper]{amsart}
\title[ Normalisers of congruence subgroups]
{  The structure of the 
 Normalisers of the congruence subgroups
 of the
 Hecke group $G_5$ 
}

\begin{document}

\baselineskip=11pt

\keywords{ congruence subgroups, Hecke groups}
\subjclass[2000]{11F06}

\maketitle

\vspace{-0.3in}
{\small $$\mbox {Mong Lung Lang}$$}

\vspace{-0.2in}

\begin{abstract} 
Let $\lambda  = 2$cos$(\pi /5)$ and let 
 $\tau \in \Bbb Z[\lambda]$. Denote 
 the normaliser of $G_0(\tau)$ of the 
 Hecke group $G_5$ in $PSL_2(\Bbb R)$
 by 
 $N(G_0(\tau))$. We prove  that
$N(G_0(\tau)) =  G_0(\tau/h)$,
 where $h$ is the largest divisor of $4$
 such that $h^2$ divides $\tau$. Further, 
$N(G_0(\tau))/G_0(\tau)$ is either 1 (if $h=1$), 
$\Bbb Z_2 \times \Bbb Z_2$ $($if $h=2)$ or
$\Bbb Z_4 \times \Bbb Z_4$ $($if $h=4)$.
\end{abstract}

\section{Introduction}
\noindent
In [LT], it is shown that $G_0(\tau) \subseteq G_5$
  is self-normalised
 if $\tau$ is square free 
 $(G_5$ is the Hecke group associated to $\lambda
 =$ 2cos$\,(\pi/5$), an element
 of $G_5$  is in $G_0(\tau)$ if and only if
 its (2,1)-entry is a multiple of $\tau$).
In this article, we complete  our study
 of the normaliser 
 of $G_0(\tau)$.
This is the first step towards the determination
 of the normalisers of the congruence 
 subgroups of all Hecke groups $G_q$.
 Note that the normalisers of the congruence
 subgroups of $G_3 = PSL_2(\Bbb Z),\, G_4,\,G_6$  have 
 been determined by Akbas-Singerman
 [AS],  Atkin-Lehner [AL], Conway [C],
 and Lang [La]. Our main result is (Section 6) :

\medskip
\noindent {\bf The Main Theorem.} {\em 
 Let $(\tau) $ be a nontrivial ideal of $\Bbb Z[\lambda]$
 and let $N(G_0(\tau))$ be the normaliser of 
 $G_0(\tau)$ in $PSL_2(\Bbb R)$. Then 
$N(G_0(\tau)) =  G_0(\tau/h)$,
 where $h$ is the largest divisor of $h_5 = 4$
 such that $h^2$ divides $\tau$. Further, 
$N(G_0(\tau))/G_0(\tau)$ is either $1$ $($if $h=1)$,
$\Bbb Z_2 \times \Bbb Z_2$ $($if $h=2)$ or
$\Bbb Z_4 \times \Bbb Z_4$ $($if $h=4)$.
}

\medskip
\noindent It is worthwhile to point out that
as an application of our main theorem,  the 
 number 4 in the main theorem enjoys the following
 property (see Appendix A) :

\smallskip
\noindent {\bf Proposition A2.} {\em
If 
 $r \in \Bbb Z[\lambda]$ is $\,G_5$-elementary
 $($see the following definition$)$, then
$r$ is a divisor of $ h_5 = 4$.}

\noindent
 {\bf Definition.} 
 $x/y \in \Bbb Q\,[\lambda]$ is in {\bf reduced form} if
{\small $\left (\begin{array}{c}
x  \\
y\\
\end{array}
\right ) $} is a column vector of some $\sigma \in G_5$.
  An element $r\in \Bbb Z[\lambda]$ 
 is called {\bf $G_5$-elementary} if $r$ satisfies
 the following
 property : Whenever $x$ and $y$ are chosen 
 such that $x/ry$ is in reduced form, 
 $x$ must admit the property that $x^2\equiv 1$
 (mod $r$).

\medskip

\noindent
 It is known that the $h$
  satisfies the above proposition for $G_6=PSL_2(\Bbb Z)$ 
 is $h_6 = 24$ (Conway [C]) and  that the number 24 plays a very
important role in the description of the normalisers
 of congruence subgroups of $G_3$, $G_4$ and $G_6$
([C], [CN], [La]). We suspect that the $h_q$ for the ring 
associated to the set of finite cusps of $G_q$ plays
some role in the determination of the normalisers and
 hope that such possibility will be considered in the 
 future.
We do not know whether Proposition A2 can be proved without
the usage of the main theorem.
  A weaker version of this proposition 
  and its direct
 proof can be found in Appendix A (Proposition A4).

\smallskip
\noindent The remaining of this article is organised
 as follows. In section 3, we give a brief introduction
 of the  reduced forms which a great detail can be found in
 [R]. We actually determine the reduced forms of certain
 members in $\Bbb Q\,[\lambda]$ which are useful
 in the determination of the normalisers of $G_0(\tau)$.
Sections 4 and 5 give supergroups  and subgroups
 of $N(G_0(\tau))$. Section 6 is devoted to the proof of the 
 main theorem.

\section {Known results}

\noindent {\bf Theorem 2.1.} (Theorem 8 of [LT])  {\em
 Let $(\tau)$ be a nontrivial  ideal of 
 $\Bbb Z[\lambda]$
 and let  $(\tau_0)$ be the maximal ideal of 
 $ \Bbb  Z [\lambda]$  such that
 $( \tau_0^2) \subseteq (\tau)$.
 Then  $N(G_0(\tau)) \subseteq G_0(\tau_0)$.
}

\medskip
\noindent  Note that if $\tau = \prod p_i^{m_i} $
where the $p_i$'s are primes, 
 then $\tau_0 = \prod p_i^{m_i -\lfloor m_i/2\rfloor      },$
where $\lfloor x\rfloor $ 
is the largest integer less than or equal to $x$.
 Throughout this article, we let 
$$[\tau] = \prod p_i^{m_i -\lfloor m_i/2\rfloor }.$$
 It is clear that if $a$ and $b$ are coprime to each other,
 then
$[ab]=[a][b]$.

\medskip
\noindent {\bf Lemma 2.2.} (Lemma 1 of [CLLT])
{\em Let $A$ be a non-zero ideal of $\Bbb Z[\lambda]$. Then
$$[G_5 : G_0(A) ] = 
 N(A)\prod_{P|A}(1 + N(P)^{-1}),$$
 where the product is over the set of all
 prime ideals $P$ which divide $A$ and $N(A)$
 is the absolute norm of $A$.}

\section {Reduced forms}
\noindent
Let $x, y \in \Bbb Z[\lambda]$. 
 $x/y$ is in {\bf reduced form} if
{\small $\left (\begin{array}{c}
x  \\
y\\
\end{array}
\right ) $} is a column vector of some $\sigma \in G_5$.
 The purpose this section is to  determine the reduced
 forms of $3/n\lambda$,  $9/n\lambda$ and $5/n\lambda$
 which will be used in the proof of Propositions 4.3 and 
 4.4.

\smallskip
\noindent
 For any $a,b \in \Bbb Z[\lambda]$ such that 
the greatest common divisor of $a$ and $b$  is a unit,
applying results of Leutbecher ([L1],[L2]),
 there exists a unique rational integer
$e(a/b)$ 
such that $$ a\lambda ^{e(a/b)}/b\lambda ^{e(a/b)}$$
is in reduced form. We shall called $e(a/b)$ 
the  {\bf reduced factor}   of $a/b$.
Since the reduced forms play a central role
 in the determination of the normaliser (see Lemma 4.1),
 we shall now give an algorithm that enables us
to determine the  reduced factor.
\smallskip

\noindent 
Let $a, b \in \Bbb Z[\lambda] \setminus \{ 0\}$
be given such that 
the greatest common divisor of $a$ and $b$  is a unit.
 Then there exists
 a unique rational integer $q$ such that

\begin{enumerate}
\item[(i)] $a = (q \lambda )b + r$,

\item[(ii)]$  -|b\lambda|/2 <|r| \le |b\lambda|/2$.
\end{enumerate}
We call such a division algorithm {\bf pseudo-Euclidean}
 (see [R] for more details). In terms of 
 matrices, the above can be written as 
 $$ 
\left (
\begin{array}{rr}
1 & -q\lambda \\
0 & 1 \\
\end{array}
\right )  
\left (
\begin{array}{r}
a  \\
b \\
\end{array}
\right )  = 
\left (
\begin{array}{r}
r  \\
b \\
\end{array}
\right )\,.$$ 
Note that 
$\left (
\begin{array}{rr}
1 & -q\lambda \\
0 & 1 \\
\end{array}
\right ) \in G$.
 Applying the  pseudo
 Euclidean algorithm repeatedly, one has,
$$a = (q_1 \lambda) b + r_1\,,$$
$$b = (q_2 \lambda) r_1 + r_2\,,$$
$$r_1 = (q_2 \lambda) r_2 + r_3\,,$$
$$..........................$$
$$..........................$$
$$..........................$$
$$r_k = (q_{k+1} \lambda) r_{k+1} + r_{k+2}\,,$$
$$r_{k+1} = (q_{k+2} \lambda) r_{k+2} + 0\,.$$
The finiteness of the algorithm is 
 governed by the fact that 
the set of cusps of $G$ is $\Bbb Q\,[\lambda]
 \cup \{ \infty \}$.
Note that in terms of matrices, the above can be written as
$$
\left (
\begin{array}{c}
a  \\
b \\
\end{array}
\right )  = A
\left (
\begin{array}{c}
r_{k+2}  \\
0 \\
\end{array}
\right )\,, \eqno (3.1)$$ 
 where $A \in G_5$.
It is clear that
$\mbox{ gcd}(a,b) =\mbox{ gcd}(b, r_1) 
 = \cdots = \mbox{ gcd}
(r_{k+1}, r_{k+2}) =  r_{k+2}$
 is a unit ($a$ and $b$ are coprime).
As 
\begin{enumerate}
\item[(i)] $|r_{k+2}| <1$,
 \item[(ii)] $\lambda$ is a primitive unit, 
\end{enumerate}
 there exists $e(a/b) \in \Bbb N \cup \{0\}$ such that 
 $|r_{k+2}| = \lambda^{-e(a/b)}\,.$
Multiplying (3.1)  by $\lambda ^{e(a/b)}$,
 one has 
$$
\left (
\begin{array}{c}
a\lambda ^{e(a/b))}  \\
b\lambda ^{e(a/b)} \\
\end{array}
\right )  = A
\left (
\begin{array}{c}
\pm1  \\
0 \\
\end{array}
\right )\,.$$ 
Since $A \in G_5$
 and $\pm 1/0$ is a reduced form,
 we conclude that $
a\lambda ^{e(a/b)}/b\lambda ^{e(a/b)} $ is the reduced form of
 $a/b$ and $e(a/b)$ is the reduced factor of $a/b$.
 Note that it is clear from the above argument 
 that 
$$e(a/b) = e(r_1/b) = \cdots  = e(r_{k+1}/r_{k+2}).\eqno(3.2)$$

\medskip
\noindent The following lemma shows that if 
 $m -  n = kr,\, r \in \Bbb Z$,  then
$e(k/n\lambda) = e(k/m\lambda)$.
 Note that this equality does not hold
 for  $e(k/n)$ and $e(k/m)$.

\medskip

\noindent {\bf Lemma 3.1.} {\em 
Let $x, \,   u,\, v \in \Bbb Z[\lambda]^{\times}$. Suppose that 
 gcd$\,(x, u)$ is a unit and that 
 $(u-v)\lambda  = xm\lambda$ for some $m \in \Bbb Z$. 
Then 
$e(x/u \lambda ) = e(x/v \lambda )$. }
\medskip

\noindent {\em Proof.}
 Let $$u \lambda = x (q_1\lambda) + r_1\,,$$
$$v \lambda = x (q_2\lambda) + r_2\,,$$
 where $q_1, \, q_2 \in \Bbb Z$, 
$ -|x|\lambda /2 < r_1,\, r_2 \le |x|\lambda /2$.
 Since  $(u-v)\lambda  = xm\lambda$, we have  $r_1 = r_2$.
 Applying equation (3.2), we have
$e(x/u\lambda ) =e(r_1/x ) =e(r_2/x ) =
 e(x/v \lambda )$.
\qed

\medskip
\noindent {\bf Example 3.2.} Applying 
 Lemma 3.1, the  following gives the  reduced 
 factors of $p^a/n\lambda$, where $p \le 11$ is a prime.
These reduced forms  will be useful in our study of 
 Propositions 4.3 and 4.4.

$$\begin{array}{llrr}
p^a/n\lambda  & & e(p^a/n\lambda) & p^a\lambda ^{ e(p^a/n\lambda)}\\
\\
 2/n \lambda  & \mbox{gcd}\,(n, 2) =1    & 2 &  2(\lambda +1)\\
 4/n \lambda  &  \mbox{gcd}\,(n, 2) =1   & 2 &  4(\lambda +1)\\
 3/n \lambda  &  \mbox{gcd}\,(n, 3) =1    & 3 &  3(2\lambda +1)\\
 9/n \lambda  &n \equiv 1, 8\,\,(\mbox{mod } 9)& 3 &  9(2\lambda +1)\\
 9/n \lambda  &n \equiv 2,4,5,7\,(\mbox{mod }  9)& 9 &  9(34\lambda +21)\\
 5/n \lambda  &  \mbox{gcd}\,(n, 5) =1    & 6 &  5(8\lambda +5)\\
 25/n \lambda  &n \equiv 1,2, 23,24\,\,(\mbox{mod }  25)& 6
 &  25(8\lambda +5)\\
 25/n \lambda  &n \equiv 3,  6,19 ,22\,\,(\mbox{mod }  25)& 
12 &  25(144\lambda +89)\\
 7/n \lambda  &\mbox{gcd}\,(n, 7) =1 & 6 &  7(8\lambda +5)\\
 11/n \lambda  &n \equiv 1,10\,\,(\mbox{mod } 11)& 6 &  11(8\lambda +5)\\
\end{array}$$

\medskip

\noindent {\em Proof.} Suppose that gcd$\,(n,3)=1$.
 It follows that  $n = 3k \pm1$.  We have
$$n\lambda = 3\cdot k\lambda \pm \lambda\,,$$
$$3 = \lambda \cdot \lambda + ( 2-\lambda)\,,$$
$$\lambda =  ( 2-\lambda)\cdot 3 \lambda + (2\lambda -3)\,,$$
$$2- \lambda =  (2\lambda -3)\cdot \lambda +0\,.$$
As  $2\lambda -3 = \lambda ^{-3}$,
 the reduced factor is $3$ and 
 the reduced form of 
 $3/n \lambda $ is 
$3\lambda ^3/n\lambda^4$.
 The rest can be calculated
 similarly. \qed

\medskip
\noindent {\bf Example 3.3.}
$e((2\lambda -1)/12) =e((2\lambda -1)/96) = 6$.
$e((2\lambda -1)/192) = 18$.

\medskip

\noindent By Example 3.3,  the reduced forms of 
$(2\lambda -1)/12$, $ (2\lambda -1)/96$ and 
$(2\lambda -1)/192$ are 
$(18\lambda +11)/12(8\lambda+5)$,
$(18\lambda +11)/96(8\lambda+5)$ and
$(5778\lambda +3571)/192(2584\lambda+1597)$ 
 respectively. The above reduced forms  will be used in the 
 remark of Lemma 5.1.

\section {Supergroups of $N(G_0(\tau))$ }

\noindent The main purpose of this section is to prove
 that $N(G_0(\tau)) \subseteq G_0(\tau/h)$, where $h$ is the
 largest divisor of 4 such that $h^2|\tau$ (Theorem 4.5).
 In order to achieve this, we need two technical lemmas
 which can be found in subsection 4.1. The main result
 (Theorem 4.5)
 is in fact a corollary of Propositions 4.3 and 4.4.
\smallskip

\subsection{Technical Lemmas.} The following two lemmas
 are useful towards the determination of a supergroup
 of $N(G_0(\tau))$.

\smallskip
\noindent {\bf Lemma 4.1.}
 {\em
 Let $\tau, [\tau ] \in \Bbb Z[\lambda]$ 
 be given as in Theorem $2.1$ and let 
 $u,\, w$ be chosen such that 
 $u/w\tau$
 is  in  reduced form. Suppose that 
gcd$\,(\tau/[\tau ],  
u^2 -1 ) = x$. 
 Then 
 $N(G_0(\tau)) \subseteq  G_0(\tau/x)$.}

\medskip

\noindent{\em Proof.}
 Since $u/w\tau$ is a reduced form,
 by results of Leutbecher ([L1], [L2]),
 $G_0(\tau) $ has an element of the form
 $$\sigma = \left (
\begin{array}{ll}
u &  v\\
w\tau   & z \\
\end{array}
\right )\,.$$
Let $ A  =
\left (
\begin{array}{ll}
a &  b\\
c & d \\
\end{array}
\right ) \in N(G_0(\tau)) \,.$
 We need to show that $A\in G_0(\tau/x)$.
 Equivalently, we need to show that 
$c$ is a multiple of $\tau/x$.
 By Theorem 2.1, $c^2$ is a multiple of $\tau$.
Since 
$$A \sigma A^{-1} =
\left (
\begin{array}{cc}
* &  *\\
cdu + wd^2 \tau -c^2v - cdz  & * \\
\end{array}
\right ) \in G_0(\tau),$$
 we conclude that 
$cdu + wd^2 \tau -c^2v - cdz
=cd(u-z) +wd^2\tau -c^2v
$ is a multiple of $\tau$.
 It is clear that the last two terms in the 
 expression are multiples of $\tau$
($c^2$ is a multiple of $\tau$).
Hence $$\tau|cd(u-z).$$ Since
 gcd$\,(c ,d) =1$ and $[\tau]|c$, we have
 gcd$\,(\tau ,d) =1$. Hence  the above is the same as 
$\tau |c(u-z).$ Since both $\tau$ and $c$ are multiples
 of $[\tau]$, we have
 $$ 
\frac{\tau}{[\tau ]} \bigg |\frac{c}{[\tau ]}(u-z).
\eqno(4.1)$$
Hence 
$$c/[\tau ]  \mbox{  is a multiple of }
\frac{\tau/[\tau ]}{ gcd\,(\tau/[\tau ], u-z)}\,. \eqno(4.2)$$

\smallskip
\noindent We shall now determine $gcd\,(\tau/[\tau ], u-z)$.
Note first that gcd$\,(u, \tau) =1$. It follows that 
$$gcd\,(\tau/[\tau ], u-z) = gcd\,(\tau/[\tau ], u(u-z)).
 $$
Since $uz-vu\tau=1$ (this is the determinant of
 $\sigma$), the right hand side of the above equation becomes
$$ gcd\,(\tau/[\tau ], u(u-z))
 = gcd\,(\tau/[\tau ], u^2 -1-vw\tau)
 = gcd\,(\tau/[\tau ], u^2 -1)
  = x. $$
Statement (4.2) now can be rephrased as 
$$c/[\tau ] \mbox{ is a multiple of } 
\frac{\tau/[\tau ]}{ x}\,.$$
 Therefore, $c$ is a multiple of $\tau/x$.
 This completes the proof of the lemma.
 \qed

\smallskip
\noindent 
{\bf Definition.} ({\bf Principal Congruence Subgroups}).
 Let
$$G(\tau)
 =\left \{  \left (
\begin{array}{ll}
a &  b\\
c  & d \\
\end{array}
\right ) \in G_5
 \,:\, 
\left (
\begin{array}{ll}
a &  b\\
c  & d \\
\end{array}
\right )\equiv
\pm I
\,\,(mod\,\,\tau)
\right \}.$$
Then $G(\tau)$ is a normal subgroup of $G_5$.
 $G(\tau)$ is called the principal 
 congruence subgroup of $G_5$.

\medskip

\noindent {\bf Lemma 4.2.} {\em  Let 
 $V = \left <G_0(3^m\nu), G(\nu) \right >$,
 where gcd$\,(3,\nu)=1$.
 Then $ V = G_0(\nu)$.
}
\medskip

\noindent {\em Proof.}
Since there
 is nothing to prove if $m = 0$,  we shall
 assume that $m \ge 1$.
 By Lemma 2.2, we have  
 $[G_0(\nu ) : G_0(3^m\nu)] = 3^{2m-2}10$.
 In order to show that $V = G_0(\nu)$, it suffices
 to show that $V$ possesses more than  $[G_0(\nu )
 : G_0(3^m\nu)]/2 = 
3^{2m-2}5$
 left $G_0(3^m\nu)$-cosets.
Let
$$\Omega = \left \{ A_{xy} = 
\left ( \begin{array}{cc}
1 & x\lambda \\
0 & 1 \\
\end{array}
\right ) 
\left ( \begin{array}{cc}
1 & 0 \\
yN(\alpha) \lambda  & 1 \\
\end{array}
\right ) \,:\, 0 \le x,\,y \le 3^m-1,\,gcd\,(3,y)=1\right \}
 .$$
 It is clear that $\Omega \subseteq V$
 and that $V$ possesses $2\cdot 3^{2m-1}$
 elements ($2\cdot 3^{m-1}$ choices of $y$
 and $3^m$ choices of $x$).
Applying the fact that $\lambda^2= \lambda+1$, one can show
 by direct  calculation 
 that 
 $$A_{xy} G_0(3^m\nu)
 = A_{x'y'} G_0(3^m\nu) \mbox{ iff }
 y \equiv  y', \,\,y^2(x-x') \equiv 0\,\, (mod \,\,3^m)
.$$
Consequently,
 members in $\Omega$
 give different $G_0(3^m\nu)$-cosets.
 Since $|\Omega| = 2\cdot 3^{2m-1}$, we conclude that 
 $V$ has at least  $2\cdot 3^{2m-1} > 3^{2m-2} 5$ left 
 $G_0(3^m\nu)$-cosets.
 This completes the proof of the lemma. \qed

\subsection {Supergroups of $N(G_0(\tau))$}
 In this section, 
 we apply our results in section 4.1 and our knowledge
 about the reduced forms of $3/n\lambda$, 
  $9/n\lambda$ and $5/n\lambda$ to give 
 some  more  accurate description of the $x$'s in Lemma 4.1.
 This is done in Proposition 4.3 and 4.4. 
 As a corollary of these two propositions, we can show
 that $G_0(\tau/h)$ (see Theorem 4.5) is a supergroup
 for $N(G_0(\tau))$.

\smallskip
\noindent {\bf Proposition 4.3.} {\em   Let $\tau = 3^m\nu$
 $(3$ is a prime in $\Bbb Z[\lambda])$,
 where   gcd$\,(3,
\nu)=1$. Then
 $N(G_0(\tau)) \subseteq N(G_0(\nu))
\subseteq G_0(\nu/\gcd(\nu/[\nu], 4)
$.
}

\smallskip
\noindent {\em Proof.} 
 We shall first prove that $N(G_0(\tau)) \subseteq N(G_0(\nu))$
: 
Let $x \in N(G_0(3^m\nu))$. We need show show that 
  $x \in N(G_0(\nu))$.
By Theorem 2.1, $x \in G_0([3^m\nu]) \subseteq G_5$.
 Since the principal
 congruence subgroup $G(\nu)$ is normal in $G_5$,
 we have $x G (\nu) x^{-1} = G(\nu)$. It follows that 
$$V =  \left <G_0(3^m\nu), G(\nu) \right >
 \subseteq
 xG_0(\nu)x^{-1} \cap G_0(\nu) \subseteq G_0(\nu). \eqno(4.3)$$
By Lemma 4.2, $V = G_0(\nu)$. 
 Expression (4.3) now implies that 
  $xG_0(\alpha)x^{-1} \cap G_0(\alpha)
 = G_0(\nu)$.
This implies that $x \in N(G_0(\nu))$.
 Hence  $$N(G_0(\tau)) \subseteq N(G_0(\nu)).$$
 We shall now prove that  $N(G_0(\nu))
\subseteq G_0(\nu/\gcd(\nu/[\nu], 4)
$ :
Let $n$ be the smallest positive rational integer
 in $(\nu)$. Then gcd$\,(n, 3) =1$. By Example
 3.2, the reduced form of 
 $3/n\lambda$  is $3\lambda^3/n\lambda^4$.
By Lemma 4.1, we have
 $$N(G_0(\nu)) \subseteq G_0(\nu/x)\,,
 \mbox{ where } x = gcd\,(\nu/[\nu], (3\lambda^3)^2-1).
 \eqno(4.4)$$
By Example 3.2, the reduced
 form of  $9/n\lambda$ is either 
$9\lambda^3/n\lambda^4$ (if $n\equiv 1,8$ (mod 9)
 or $9\lambda^9/n\lambda^{10}$ (if $n\equiv 2,4,5,7$ (mod 9).
\smallskip

\noindent 
(i) Suppose that  $n\equiv 1,8$ (mod 9). By Lemma 4.1,
 we have
 $N(G_0(\nu)) \subseteq G_0(\nu/y)$,
 where $$ y= gcd\,(\nu/[\nu], (9\lambda^3)^2-1).$$
Applying (4.4), we have
 $N(G_0(\nu)) \subseteq G_0(\nu/y) \cap  G_0(\nu/x).$
 Since $G_0(u)\cap G_0(v) = G_0(\mbox{lcm}\,(u,v))$ and 
 gcd$\,(  (3\lambda^3)^2-1 ,  (9\lambda^3)^2-1     ) =4$,
 we have 
 $$N(G_0(\nu))\subseteq
G_0(\nu/x) \cap G_0(\nu/y) =
  G_0(\nu/ gcd\,(\nu/[\nu], 4)).$$

\smallskip
\noindent 
(ii) Suppose that  $n\equiv 2,4,5,7$ (mod 9). By Lemma 4.1,
 we have
$N(G_0(\nu)) \subseteq G_0(\nu/z)$,
 where $$  z = gcd\,(\nu/[\nu], (9\lambda^9)^2-1).$$
Applying (4.4), we have
 $N(G_0(\nu)) \subseteq G_0(\nu/z) \cap  G_0(\nu/x).$
 Since $G_0(u)\cap G_0(v) = G_0(\mbox{lcm}\,(u,v))$ and 
 gcd$\,(  (3\lambda^3)^2-1 ,  (9\lambda^9)^2-1     ) =4$,
 we have
 $$N(G_0(\nu)) \subseteq 
G_0(\nu/z) \cap  G_0(\nu/x)=
G_0(\nu/ gcd\,(\nu/[\nu], 4)).
\eqno\qed$$

\smallskip
\noindent 
Similar to Proposition 4.3, we have the following
 proposition for the decomposition
$\tau = \sqrt 5^m\nu$
(note that  the reduced form of 
 $5/n\lambda$ is $5\lambda^6/n\lambda^7$
 (Example 3.2) and that  
 $(5\lambda^6)^2-1 = 16(225\lambda +139)$, where
the absolute norm of  $225\lambda +139 $ is 29).

\smallskip
\noindent {\bf Proposition 4.4.} {\em 
 Let $\tau = {\sqrt5}\,\,
^m\nu$,
 where   gcd$\,(\sqrt 5,
\nu)=1 $ $($$5$ is not a  prime$)$. Then
 $N(G_0(\tau)) \subseteq N(G_0(\nu))
\subseteq G_0(\nu/\gcd(\nu/[\nu], 16(225\lambda +139))
$.}

\smallskip
\noindent 
{\bf Theorem 4.5.} 
{\em Let $(\tau) $ be a nontrivial
 ideal of $\Bbb Z[\lambda]$ and let $N(G_0(\tau))$ be the 
 normaliser of $G_0(\tau)$ in $PSL_2(\Bbb R)$. Then
 $N(G_0(\tau)) \subseteq G_0(\tau/h)$,
 where $h$ is the largest divisor of $4$
 such that $h^2$ is a divisor of $\tau$.}

\smallskip
\noindent
{\em Proof.} 
 Let  $\tau = 2^a3^b\sqrt 5 \,^c \alpha$, where
 gcd$\,(\alpha ,  6\sqrt 5)=1$.
Applying Propositions
 4.3 and 4.4,
 we have 
$$N(G_0(\tau)\subseteq
 G_0\left (\frac{2^a\sqrt 5 \,^c \alpha}{gcd\,\left (
\frac{2^a\sqrt 5 \,^c \alpha}{[2^a\sqrt 5\,^c \alpha  ]}
, 4 \right )}\right )
\cap  G_0\left (\frac{2^a 3^b \alpha}{gcd\,\left (
\frac{2^a 3^b\alpha }{[2^a 3^b\alpha]}
, 16(225\lambda +139 )\right ) } \right ).
 $$
Since $[xy] = [x][y]$ whenever gcd$\,(x,y)=1$,
$ 225\lambda +139$ has absolute norm 29
and $G_0(x) \cap G_0(y) =G_0(\mbox{lcm}\,(x,y))$, the
  above 
 can be simplified into 
$$ N(G_0(\tau)) \subseteq 
G_0\left (\frac{2^a 3^b\sqrt 5\, ^c \alpha}{gcd\,(2^{a-\lfloor a
\rfloor}, 4 )}\right )
 =  G_0\left (\frac{\tau}{gcd\,(2^{a-\lfloor a \rfloor},
4)}\right ).\eqno (4.5)$$
Note that  $gcd\,(2^{a-\lfloor a \rfloor}, 4)$ 
in equation (4.5)
 is nothing but the largest divisor $h$ of 4 such that 
 $h^2|\tau$. This completes the proof of the theorem. \qed

\section{ Subgroups of $N(G_0(\tau))$
}

\noindent The  aim of this section is to 
 find  large subgroups of $N(G_0(\tau))$.
 It is proved in Lemma 5.2 that the supergroups
 for $N(G_0(\tau))$ (Theorem 4.5) are actually
 the subgroups of  $N(G_0(\tau))$ as well.

\medskip
\noindent 
Recall that $\Bbb Z[\lambda]$ is a principal
 ideal domain and 
 $2$ is a prime in $\Bbb Z[\lambda]$.
 An element $\tau$  in $\Bbb Z[\lambda]$
 is  even if $\tau = 2 x$
 for some $x \in \Bbb Z[\lambda]$.
Unlike the rational integers $\Bbb Z$, 
 the sum of two odd integers in  $\Bbb Z[\lambda]$
 does not have to be even (both 
 1 and $\lambda $ are odd, yet the sum of
 $1$ and $\lambda$ is not even).
As a consequence, the trace of 
  $\sigma \in PSL_2(\Bbb Z[\lambda])$,
where the 
(2,1) entry of $\sigma$ is even,
 does not have to be even
 (see remark of lemma 5.1). However, as the following 
 lemma shows, the trace of $\sigma \in G_5$ is even
 if the (2,1) entry of 
 $\sigma$ is even.

\medskip

\noindent {\bf Lemma 5.1.} {\em 
 Let $\left (
\begin{array}{rr}
a  & b \\
2c & d \\
\end{array}
\right ) \in G_0(2)$. Then 
 $a + d$,  $a-d$ 
are even and  $a^2-1$ is a multiple of $4$.
}

\medskip

\noindent {\em Proof.} An easy study of a fundamental
 domain of $G_0(2)$ shows that 
 $G_0(2)$ can be  generated by
$$ 
\left (
\begin{array}{rr}
1  & \lambda \\
0 & 1 \\
\end{array}
\right ) \,,
\left (
\begin{array}{rr}
2 \lambda+1  & - \lambda-2 \\
2 \lambda+2 & -2 \lambda-1 \\
\end{array}
\right ) \,,
\left (
\begin{array}{rr}
2 \lambda+1  & - \lambda \\
2 \lambda & -1 \\
\end{array}
\right ) $$
 (see [LLT]).
Note that 

\begin{enumerate}
\item[(a)] the sum and 
difference between the (1,1) and (2,2)
 entries 
of the matrices  above 
are even, (b)  $(2\lambda +1)^2 -1$
 is a multiple of 4.
\end{enumerate}
 As a consequence, the generators of $G_0(2)$ satisfy
 the conclusion of our lemma. Since 
 every element
 of $G_0(2)$ can be written as a word of  the generators
 of length $n$ for some $n$,
  the lemma   can be proved
 by applying induction on $n$.
\qed

\medskip

\noindent {\bf Remark.}  (a)
  An easy study of  
{\small $ 
\left (
\begin{array}{rr}
3 \lambda-1  & \lambda \\
2 \lambda &  \lambda \\
\end{array}
\right ) 
$}
shows that 
Lemma 5.1 does not hold for $PSL_2(\Bbb Z [\lambda])$.
\medskip

\noindent (b) The squares of the (1,1) entries 
 of the following matrices  (see Example 3.3) are not congruent
 to 1 modulo 8.
 This means that  
 the second part of Lemma 5.1 ($a^2 -1\equiv 0$
 (mod 4)) can not be improved for the following groups.
   $$ 
\left (
\begin{array}{rr}
6 \lambda+5  & -\lambda  \\
6\lambda &  -1 \\
\end{array}
\right ) \in G_0(6)\,,\,\,
\left (
\begin{array}{rr}
18\lambda+11  & *  \\
12(8\lambda+5) &  * \\
\end{array}
\right ) \in G_0(12)\,,$$

$$\left (
\begin{array}{rr}
18\lambda+11  & *  \\
96(8\lambda+5) &  * \\
\end{array}
\right ) \in G_0(96)\,, \,\,
\left (
\begin{array}{rr}
5778\lambda+3571  & *  \\
192(2584\lambda+1597) &  * \\
\end{array}
\right ) \in G_0(192)\,.$$

\medskip

\noindent {\bf Lemma  5.2.} {\em
 Let $\tau \in \Bbb Z[\lambda]$. Suppose that 
 ${\tau}$ is even. Then the following holds.
\begin{enumerate}
\item[(a)]
  If $4|\tau$, then
$G_0(\tau/2
) \subseteq  N(G_0(\tau))$.
\item[(b)]
 If $16|\tau$, then
$
G_0(\tau/4
) \ \subseteq  N(G_0(\tau))$.
\end{enumerate}}

\medskip
\noindent {\em Proof.} 
 (a)
 Let 
 $
A = \left (
\begin{array}{cc}
a & b \\
c\tau/2 & d \\
\end{array}
\right ) \in G_0(\tau/2)$, 
$B =
\left ( \begin{array}{cc}
u & v \\
w\tau & x \\
\end{array}
\right ) \in G_0(\tau)$. Then

 $$  ABA^{-1} = 
\left (
\begin{array}{cc}
*  & * \\
cd(u-x)\tau/2 -c^2v\tau^2/4 +wd^2\tau  & * \\
\end{array}
\right )
 \,.$$
 Since  $u-x$ is even (Lemma 5.1),
 ${\tau}$ is 
a multiple of 4 (our assumption), we conclude that
 the three terms in  
 $cd(u-x)\tau/2 -c^2v\tau^2/4 +wd^2\tau $
 are multiples of $\tau$. Hence
 $cd(u-x)\tau/2 -c^2v\tau^2/4 +wd^2\tau $
 is a multiple of $\tau$.
 As a consequence, $ABA^{-1} \in G_0(\tau)$.
 This implies that 
$
G_0(\tau/2
) \subseteq N(G_0(\tau))\,.$

\medskip
\noindent (b) 
Let 
 $
A = \left (
\begin{array}{cc}
a & b \\
c\tau/4 & d \\
\end{array}
\right ) \in G_0(\tau/4)$, 
$B =
\left ( \begin{array}{cc}
u & v \\
w\tau & x \\
\end{array}
\right ) \in G_0(\tau)$. Then 
 $$  ABA^{-1} = 
\left (
\begin{array}{cc}
*  & * \\
cd(u-x)\tau/4 -c^2v\tau^2/16 +wd^2\tau  & * \\
\end{array}
\right )
 \,.$$
Since 
\begin{enumerate}
\item[(i)]  $u -x = (u^2-1 - vw\tau)/u$
   (det$\,B  = ux-wv\tau=1),$
\item[(ii)] $u^2-1$ is a multiple of 4 (Lemma 5.1),
\item[(iii)] $\tau$ is a multiple of 16
 (our assumption),
\end{enumerate}
 we conclude that 
$u-x$ is a multiple of 4. 
Hence all three terms in $cd(u-x)\tau/4 -c^2v
\tau^2/16 +wd^2\tau $
 are multiples of $\tau$.
 This implies that 
$cd(u-x)\tau/4 -c^2v\tau^2/16 +wd^2\tau $
 is a multiple of $\tau$.
As a consequence, $ABA^{-1} \in G_0(\tau)$.
 Thus
$
G_0(\tau/4
) \subseteq N(G_0(\tau))\,.$
\qed

\section{ The Main Theorem : Normaliser of $G_0(\tau)$}

\noindent The following theorem
 follows immediately by
 applying  Theorem 4.5 and Lemma 5.2.
\smallskip

\noindent {\bf Theorem 6.1.} {\em 
Let $(\tau) $ be a nontrivial
 ideal of $\Bbb Z[\lambda]$ and let $N(G_0(\tau))$ be the 
 normaliser of $G_0(\tau)$ in $PSL_2(\Bbb R)$. Then
 $N(G_0(\tau)) = G_0(\tau/h)$,
 where $h$ is the largest divisor of $4$
 such that $h^2$ is a divisor of $\tau$.
}
\medskip

\noindent 
The following proposition determines the group structure
 of $N(G_0(\tau))/G_0(\tau)$.
\smallskip

\noindent {\bf Proposition 6.2.} {\em 
$N(G_0(\tau))/G_0(\tau) = G_0(\tau/h)/G_0(\tau) $ is either
 $1$ $($if $h=1)$,
$\Bbb Z_2 \times \Bbb Z_2$ $($if $h=2)$ or
$\Bbb Z_4 \times \Bbb Z_4$ $($if $h=4)$.}

\medskip
\noindent
{\em Proof.}  Recall first that $h^2|\tau$.
 There is nothing to prove for the case
 $h=1$. 
 In the case $ h= 2$, $[N(G_0(\tau)) : G_0(\tau)]=
[G_0(\tau/2) : G_0(\tau)]=  4$ (Lemma 2.2).
For any $\sigma \in G_0(\tau/2)$, since $\tau/2$ is even,
we have 
 $\sigma = \left (
\begin{array}{ll}
a &  b\\
2c & d \\
\end{array}
\right ) 
 \in G_0(\tau/2) \subseteq G_0(2)$. Applying
 Lemma 5.1, the (2,1)-entry of $\sigma^2$,
 which is $2c(a+d)$,  is a multiple of 
 $\tau$. Hence $\sigma ^2 \in G_0(\tau)$. In particular,
every nontrivial element in $G_0(\tau/2)/G_0(\tau)$ is 
 of order 2. Hence $N(G_0(\tau)) /G_0(\tau)=
G_0(\tau/2) /G_0(\tau) \cong \Bbb Z_2 \times \Bbb Z_2$.
We shall now assume that $h = 4$.

\medskip
\noindent (i) $[N(G_0(\tau)) : G_0(\tau)]=
[G_0(\tau/4) : G_0(\tau)]=  16$ and $G_0(\tau/4)/G_0(\tau)$ is 
 abelian.
\medskip

\noindent {\em Proof of}  (i). By Lemma 2.2,
$[G_0(\tau/4) : G_0(\tau)]=  16$. Let
$$\sigma_1 = 
 \left (
\begin{array}{cc}
a_1 & b_1 \\
c_1\tau/4 & d_1 \\
\end{array}
\right ),\,\, \sigma_2 = 
 \left (
\begin{array}{cc}
a_2 & b_2 \\
c_2\tau/4 & d_2 \\
\end{array}
\right ).$$
 The (2,1)-entry of $\sigma_1\sigma_2\sigma_1^{-1}
 \sigma_2^{-1}$ is
$$\frac{\,\,\tau\,\,}{4}[(a_2c_1+d_1c_2)(c_2b_1\tau/4 +d_1d_2)
 - (a_1c_2+d_2c_1)(c_1b_2\tau/4 +d_1d_2)]. \eqno(6.1)$$
Since $16$ divides $\tau$, 
$$(6.1) \equiv\frac{\,\,\tau\,\,}{4}
 d_1d_2((a_2-d_2)c_1 +c_2(d_1-a_1))
 \,\,\,(\mbox{mod}\,\, \tau).\eqno(6.2)$$
Applying Lemma 5.1, $a_i^2 \equiv 1$ ( mod 4).
Since $a_id_i\equiv 1$ (mod 4), we conclude that 
 $a_i -d_i \equiv 0$ (mod 4).
Hence the right hand side of (6.2) is a multiple of $\tau$.
 Hence 
  $(6.1) \equiv 0$ (mod $\tau$).
 Equivalently, the (2,1)-entry
of  $\sigma_1\sigma_2\sigma_1^{-1}\sigma_2^{-1}
  $ is a multiple
 of $\tau$. Therefore, 
 $\sigma_1\sigma_2\sigma_1^{-1}
 \sigma_2^{-1} \in G_0(\tau)$ for all $\sigma_1, \sigma_2 \in 
 G_0(\tau/4)$. Hence 
 $G_0(\tau/4)/G_0(\tau)$ is abelian.
\medskip

\noindent (ii) $G_0(\tau/4)/G_0(\tau) \cong \Bbb Z_4
 \times \Bbb Z_4.$    {\em Proof of } (ii).
 Let
$$\sigma = 
 \left (
\begin{array}{cc}
a & b \\
c\tau/4 & d \\
\end{array}
\right ).$$
 Applying Lemma 5.1, the (2,1)-entry
 of $\sigma^4$ is a multiple of $\tau$.
 It follows that $\sigma^4 \in G_0(\tau)$. Hence every
 element in $G_0(\tau/4)/G_0(\tau)$ is of order 1, 2 or 4.
 Let $n/4$ be the smallest positive rational integer in $
(\tau/4)$. The order of {\small  $
 \left (
\begin{array}{cc}
1 & 0 \\
 n \lambda/4  & 1\\
\end{array}
\right )$} is 4 in  $G_0(\tau/4)/G_0(\tau)$. 
 Since $G_0(\tau/4)/G_0(\tau)$ is abelian and every
 element in $G_0(\tau/4)/G_0(\tau)$ is of order 1, 2 or 4,
 we have 
   $$G_0(\tau/4)/G_0(\tau) \cong \Bbb Z_4
 \times \Bbb Z_4\mbox{ or }G_0(\tau/4)/G_0(\tau) \cong \Bbb Z_4
 \times \Bbb Z_2 \times \Bbb Z_2.$$ 
 Suppose that we are in the latter case. It follows
 that a representative  $\sigma$ of 
 some  nonidentity element in  $G_0(\tau/4)/
G_0(\tau/2) \cong \Bbb Z_2 \times \Bbb Z_2$ is of order
 2 in $G_0(\tau/4)/G_0(\tau)$. Let 
{\small $\sigma = 
 \left (
\begin{array}{cc}
a & b \\
c\tau/4 & d \\
\end{array}
\right )$  }  be such an element.
Since $\sigma$ is of order 2 in $G_0(\tau/4)/G_0(\tau)$,
 Then the (2,1)-entry
 of $\sigma^2$ must be a multiple of $\tau$.
 Direct calculation shows that the (2,1)-entry of 
 $\sigma^2$ is $c(a+d)\tau/4$.
Since $c$ is odd (this follows
 from the fact that 
 $\sigma \notin G_0(\tau/2))$, we have $a+d$ is a multiple of 4.
  Since $ad \equiv 1$ (mod 4) (determinant of $\sigma$ is 1),
we have $a^2 \equiv -1$ (mod 4).
 A contradiction (see Lemma 5.1).
  Hence such $\sigma$ does not exist. It follows that 
 $G_0(\tau/4)/G_0(\tau) \cong \Bbb Z_4
 \times \Bbb Z_4.$\qed

\medskip

\noindent {\bf Remark.} Unlike the modular group case,
where $  \Gamma _0(n/h) \subseteq N(\Gamma_0( n)) $
 ($h$ is the largest divisor of 24 such that $h^2$ divides
 $n$), the number 3 does not play any role in the normaliser
 of $G_0(\tau) \subseteq G_5$.
 This can be seen in Theorem 6.1 or the following
 example, which shows that $N(G_0(9))$ does not contain
{\small  $ \left (
\begin{array}{cc}
1 & 0 \\
3\lambda  & 1 \\
\end{array}
\right )  $}.

\smallskip
\noindent 
 {\bf Example.} By Example 3.2, the reduced form 
 of $4/9\lambda$ is $4\lambda ^2/9\lambda ^3$.
 It follows that $G_0(9)$ has an element of the form
{\small  $ \sigma =  \left (
\begin{array}{cc}
4\lambda ^2 & a \\
 9\lambda ^3 & b \\
\end{array}
\right ). $} Direct calculation shows that 
 {\small $$  \left (
\begin{array}{cc}
1 & 0 \\
 3\lambda  & 1 \\
\end{array}
\right )  \left (
\begin{array}{cc}
4\lambda ^2 & a \\
 9\lambda ^3 & b \\
\end{array}
\right )
 \left (
\begin{array}{cc}
1 & 0 \\
- 3\lambda  & 1 \\
\end{array}
\right )
= \left (
\begin{array}{cc}
* & * \\
21 \lambda ^3 -9a\lambda ^2 -3b \lambda & * \\
\end{array}
\right ).
$$}
 
\noindent
Since the determinant of $\sigma$ is 1,
$4b\lambda ^2 \equiv 1$ (mod 9). Hence 
$21 \lambda ^3 -9a\lambda ^2 -3b \lambda$
 is not a multiple of 9. In particular,
 {\small $\left (
\begin{array}{cc}
1 & 0 \\
3\lambda  & 1 \\
\end{array}
\right )$} does not normalise   $G_0(9)$.

\medskip
\begin{center}
{\bf Appendix A}
\end{center}

\medskip
\noindent

\smallskip

\noindent
 {\bf Definition A1.} 
  An element $r\in \Bbb Z[\lambda]$ 
 is called {\bf $G_5$-elementary} if $r$ satisfies
 the following
 property : Whenever $x$ and $y$ are chosen 
 such that $x/ry$ is in reduced form, 
 $x$ must admit the property that $x^2\equiv 1$
 (mod $r$).

\smallskip
\noindent {\bf Proposition A2.} {\em
If $r \in \Bbb Z[\lambda]$ is $G_5$-elementary, then
$r$ is a divisor of $ h_5 = 4$.}

\medskip
\noindent 
 {\em Proof.} Let $r$ be elementary and let
 $
A = \left (
\begin{array}{cc}
a & b \\
cr & d \\
\end{array}
\right ) \in G_0(r)$, 
$B = \left ( \begin{array}{cc}
u & v \\
wr^2 & x \\
\end{array}
\right ) \in G_0(r^2).$ Then

 $$  ABA^{-1} = 
\left (
\begin{array}{cc}
*  & * \\
cd(u-x)r -c^2vr^2 +wd^2 r^2  & * \\
\end{array}
\right )
 \,.$$
 Since  $ux \equiv 1$ (mod $r^2$) and $r$ is elementary,
 we have $u-x \equiv 0$ (mod $r$).
  Hence  
 $cd(u-x)r -c^2vr^2 +wd^2r^2 $
 is a multiple of $r^2$.
 As a consequence, $ABA^{-1} \in G_0(r^2)$.
 This implies that $G_0(r)= G_0(r^2/r) 
 \subseteq N(G_0(r^2)).$   By Theorem 6.1,
$  N(G_0(r^2)) = G_0(r^2/h)$, where $h$ is the 
 largest  divisor
 of $4$ such that $h^2|r^2$. 
 It follows that  $G_0(r) \subseteq  G_0(r^2/h)$.
This implies that $r$ is a divisor of 4. \qed

\medskip
\noindent {\bf Discussion.}
 We do not have a direct proof of Proposition A2. 
 The following weaker version about the reduced
 forms of $G_5$ does enjoy a proof without the usage
 of the main theorem.

\smallskip
\noindent {\bf Definition A3.} 
 $r \in \Bbb
 Z[\lambda]$ is called
{\bf strongly $G_5$-elementary} if every divisor
 of $r$ is $G_5$-elementary.

\smallskip
\noindent

\noindent {\bf Proposition A4.} {\em Suppose that  $r$
 is strongly $G_5$-elementary. 
 Then $r$ is a divisor of $h_5 = 4.$}

\medskip
\noindent {\em Proof.} Suppose that $r$ has a
 non-unit  odd divisor
 $r_1$.
 Let $N(r_1) = r_1s_1 = n\in \Bbb Z$.
By Example 3.2, the reduced form of 
$2 /r_1s_1$ is   $2\lambda^2/n\lambda^3$.
 By our assumption, 
\smallskip

\begin{center}
$(2\lambda^2)^2 
=(2(\lambda+1))^2  = 12\lambda +8 \equiv 1$ (mod $r_1)$.
\end{center}

\smallskip
\noindent 
Hence  $ r_1 = 12\lambda +7$    ($12\lambda +7$ is a prime
 of norm 11).
By Example 3.2, the reduced form of $3/r_1s_1\lambda$ is
  $3\lambda^3/11\lambda^4$.
Direct calculation shows that $(3\lambda^3)^2 \not\equiv 1$
 (mod  $r_1$). A contradiction. Hence $r$ must be a power
 of 2.  Let
 $r= 2^m$. Applying Lemma 5.1, $r=2^2$ admits the above
 mentioned property. 
Since  $(2\lambda+1)/2\lambda$ is in reduced
 form (see Lemma 5.1) and 
$(2\lambda+1)^2 \not\equiv 1$ (mod $2^3)$,
$8$ is not $G_5$-elementary.
 Hence  $r$ must be  a divisor of $h_5 = 4$. \qed

\bigskip

{\small

\noindent DEPARTMENT OF MATHEMATICS,\\
NATIONAL UNIVERSITY OF SINGAPORE,\\
SINGAPORE 117543,\\
REPUBLIC OF SINGAPORE\\ }

\noindent {\tt e-mail: matlml@math.nus.edu.sg}

\bigskip


\end{document}